\documentclass[12pt,a4paper]{article}
\usepackage{amssymb}
\usepackage{amsmath}

\newtheorem{theorem}{Theorem}[section]
\newtheorem{definition}[theorem]{Definition}
\newtheorem{lemma}[theorem]{Lemma}
\newtheorem{corollary}[theorem]{Corollary}

\date{}

\begin{document}

\title{Lyndon-Shirshov  basis and anti-commutative algebras\footnote{Supported by the NNSF of China (Nos. 11171118;
10911120389) and the NSF of Guangdong Province (No.
S2011010003374).}}
\author{
L. A. Bokut\footnote {Supported by RFBR 09-01-00157,
LSS--3669.2010.1 and SB RAS Integration grant No. 2009.97 (Russia)
and Federal Target Grant ¡°Scientific and educational personnel of
innovation Russia¡± for 2009-2013
(government contract No. 02.740.11.5191).} \\
{\small \ School of Mathematical Sciences, South China Normal
University}\\
{\small Guangzhou 510631, P. R. China}\\
{\small Sobolev Institute of Mathematics, Russian Academy of
Sciences}\\
{\small Siberian Branch, Novosibirsk 630090, Russia}\\
{\small Email: bokut@math.nsc.ru}\\
\\
 Yuqun
Chen\footnote {Corresponding author.} \  and Yu Li\\
{\small \ School of Mathematical Sciences, South China Normal
University}\\
{\small Guangzhou 510631, P. R. China}\\
{\small Email: yqchen@scnu.edu.cn}\\
{\small LiYu820615@126.com}}

\maketitle

\noindent\textbf{Abstract:} Chen, Fox, Lyndon 1958 \cite{CFL58} and
Shirshov 1958 \cite{Sh58} introduced non-associative Lyndon-Shirshov
words and proved that they form a linear basis of a free Lie
algebra, independently. In this paper we give another approach to
definition of Lyndon-Shirshov basis, i.e., we find an
anti-commutative Gr\"{o}bner-Shirshov basis $S$ of a free Lie
algebra such that $Irr(S)$ is the set of all non-associative
Lyndon-Shirshov words, where $Irr(S)$ is the set of all monomials of
$N(X)$, a basis of the free anti-commutative algebra on $X$,  not
containing maximal monomials of polynomials from  $S$. Following
from Shirshov's anti-commutative Gr\"{o}bner-Shirshov bases theory
\cite{S62a2}, the set $Irr(S)$ is a linear basis of a free Lie
algebra.

\noindent \textbf{Key words: } Lie algebra, anti-commutative
algebra, Lyndon-Shirshov words, Gr\"{o}bner-Shirshov basis

\noindent \textbf{AMS 2000 Subject Classification}: 16S15, 13P10,
17Bxx

\section{Introduction}

The first linear basis of a free Lie algebra $Lie(X)$ had been given
by M. Hall 1950 \cite{M.H50}. He proved that P. Hall long
commutators 1933 \cite{P.H33} form a linear basis of the algebra.
Let us remind that P. Hall (non-associative) monomials are  $x_i\in
X, x_i>x_j$ if $i>j,$ and then, by induction on degree (length), in
$((u)(v))$ both $(u)$ and $(v)$ are Hall monomials with $(u)>(v)$,
and in $(((u_1)(u_2))(v))$ it must be $(u_2)\leq (v)$. He used a
deg-ordering of  monomials - monomial of greater degree is greater.
In what follow ``monomial'' would mean ``non-associative monomial''
(``non-associative word'' in Kurosh's terminology \cite{Ku47},
adopted by Shirshov, see \cite{S09}).

A.I. Shirshov in his Candidate Science Thesis 1953 \cite{S53a},
published in 1962  \cite{S62a} (it was the first issue of a new
Malcev's Journal ``Algebra and Logic"), generalized this result and
found a series of bases of $Lie(X)$ depending on any ordering of Lie
monomials $(u)'$s on X such that $((u_1)(u_2)) > (u_2)$ (it was
rediscovered by Viennot 1973 \cite{V73}, see also 1978 \cite{V78}).
Any deg-ordering  of course satisfies the above condition. Now this
series of bases of $Lie(X)$ is called Hall-Shirshov or even Hall
bases (the later is used in \cite {Reu1993}).

In 1963 $\cite{B63}$\footnote {Actually, the first author found the
result in 1959 but could not published it before Shirshov's one. At
last Shirshov said to the first author: ``Your result pushes me to
publish my result".} the first author  found another ordering of
non-associative monomials with Shirshov's condition and as a result
he found a linear basis of $Lie(X)$ compatible with the lower
central series of $Lie(X)$ (the result was rediscovered by C.
Reutenauer, see his book \cite{Reu1993}, Ch. 5.3) and even a basis
compatible with any polynilpotent series
$$
L\supset  L^{n_1} \supset (L^{n_1})^{n_2}\supset \dots (\dots
((L^{n_1})^{n_2})\dots )^{n_k}\supset \dots,\ \ n_i\geq 2.
$$

It was a breakthrough in the subject that Shirshov 1958 and
Chen-Fox-Lyndon 1958 (in the same year!) found a new basis of
$Lie(X)$ now called Lyndon or Lyndon-Shirshov basis. Actually it is
an example of the Shirshov's series of bases relative to
lex-ordering of monomials: $(u)>_{lex}(v)$ if $u>_{lex}v$ for
corresponding (associative) words $u,v$, otherwise one needs to
compare the first monomials of $(u)$, $(v)$, then the second ones.
For example, $a>[ab]$, $[a[[ab]b]]> [[a[ab]]b]$. If one starts to
construct the Hall-Shirshov basis of $Lie(X)$ using the lex-ordering
then one will get Lyndon-Shirshov basis automatically. For example,
let $X=\{a,b\}$, $a>b$. Then the first Lyndon-Shirshov monomials are
$a,b, [ab], [[ab]b]=[abb], [a[ab]]=[aab], [a[[ab]b]]=[aabb],
[a[a[ab]]]=[aaab], [[a[ab][ab]]]=[aabab], [a[a[a[ab]]]]=[aaaab],
\dots, [[a[abb]][ab]]=[aabbab]$, \dots.

We would see that if $[u]$ is a Lyndon-Shirshov monomial then the
underlying word $u$ is a Lyndon-Shirshov word in the sense that
$u=vw> wv$ lexicographically for any non-empty words $v, w$. There
is a one to one correspondence between Lyndon-Shirshov  monomials
and Lyndon-Shirshov words. If $u$ is a Lyndon-Shirshov (LS  for
short) word then Shirshov 1953 \cite{S53b} and Lazard 1960
\cite{La1960} elimination process gives rise to a Lyndon-Shirshov
Lie monomial. For example, one has a series of bracketing starting
with LS word $aaab$: $ aa[ab]$ (join the minimal letter $b$ to the
previous one), $a[a[ab]]$ (again join the minimal new ``letter"
$[ab]$ to the previous letter  $a$), $[a[a[ab]]]$ (again join
$[a[ab]] $ to $a$). At last $[a[a[ab]]]$ is a Lyndon-Shirshov  Lie
monomial.

Lyndon-Shirshov basis became one of popular bases of free Lie
algebras (cf. for example, \cite{BK94,Reu1993}). One of the main
applications of Lyndon-Shirshov basis  is the Shirshov's theory of
Gr\"{o}bner-Shirshov bases theory for Lie algebras \cite{S62b}.

Original Shirshov 1958 \cite{Sh58} definition is as follows. Let
$X=\{x_i, i=1,2,...\}$ be a totally ordered alphabet, $x_i>x_j$ if
$i>j$. A non-empty monomial $u\in X^*$ is called regular
(``pravil'noe" in Russian) if $u=vw> wv$ for any non-empty words $v,
w$. Let $<$ be the lexicographical ordering of regular words. A
monomial $[u]$ is called regular if (1) the word $u$ is regular, (2)
if $[u]=[[v][w]]$, then $[v], [w]$ are regular words, (3) if in (2)
$[v]=[[v_1][v_2]]$, then $v_2\leq w$.

It is easy to see that regular monomials are defined inductively
starting with ${x_i}$. From condition (2) it follows that $v>w$.
From (1)-(3) it easily follows that $w$ is the longest proper
regular suffix of $u$, see \cite{Sh58a}. As we explained above, in
view of his thesis 1953, the definition of regular words and regular
monomials were very natural for Shirshov since he decided to use the
lexicographical ordering of monomials and underlying words. Actually
Shirshov understood all properties of regular words and monomials by
induction on degree using the Lazard-Shirshov elimination process,
see, for example, \cite{BC07}.

Standard  words ($u=vw< wv$ lexicographically for any non-empty
$v,w$; it is the same as regular word if we invert the ordering on
$X$) were defined by Lyndon  1954 \cite{L54}.

Original Chen-Fox-Lyndon 1958 \cite{CFL58} definition of the linear
basis is (also) based on the notion of standard word and its
bracketing $[u]=[[v][w]]$, where $w$ is the longest proper standard
suffix of $u$ (then $v$ is automatically standard).

Unfortunately the both Lyndon's and Chen-Fox-Lyndon's papers
\cite{CFL58,L54} were unknown to many authors for many years. Many
authors, started with 1958, call the basis and words following
Shirshov as regular basis and words, see, for example,
\cite{BaMPZa92,BaO11,Bo72,C65,C68,K08,Mikhalev,Ra89,Ra94,U95,Z92}.
We and some authors call them Lyndon-Shirshov words, see, for
example, \cite{BC07,BK94,Ci06}.

Many authors, started with 1983, call the words and monomials as
Lyndon words and Lyndon basis, see, for example,
\cite{Lo83,Reu1993}.

It is of some interest to cite some early papers on the matter by M.
Sch\"utzenberger.

In Sch\"utzenberger-Sherman 1963 \cite{SchS63}, both the Shirshov
1958 and the Chen-Fox-Lyndon 1958 papers are cited. What is more,
they formulated and  used a result, Lemma 2, in Shirshov 1958
\cite{Sh58}\footnote {From Sch\"utzenberger-Sherman 1963
\cite{SchS63}: ``We are indebted to P.M. Cohn for calling our
attention to the Shirshov paper. Every $f\in F^+$ has one and only
one factorization $f=h_{i_1}h_{i_2}...h_{i_m}$, where $h_{i_j}$
belong to $H$  and satisfy $h_{i_m}\geq \dots \geq h_{i_2}\geq
h_{i_1}$." Here $H$ is the set of Lyndon-Shirshov words in the
alphabet $F$.} (but they did not claim (!) it is a Chen-Fox-Lyndon
1958 paper result; actually, it is not, see below). The result is
sometimes called as ``Central result on Lyndon words" (see for
example Springer Online, Encyclopedia of Mathematics (edited by
Michiel Hazewinkel)).

In Sch\"utzenberger 1965 \cite{Sch65}, both Chen-Fox-Lyndon 1958 and
Shirshov 1958 papers are again cited in according with ``Central
result on  Lyndon words". As it follows from Berstel-Perrin's paper
\cite{BeP07}\footnote{From \cite{BeP07}: ``A famous theorem
concerning Lyndon words asserts that any word $w$ can be factorized
in a unique way as a non-increasing product of Lyndon words, i.e.,
$w=x_1x_2\dots x_n$ with $x_1\geq x_2\geq \dots \geq x_n$. This
theorem has imprecise origin. It is usually credited to
Chen-Fox-Lyndon, following the paper of Sch\"utzenberger 1965 [120]
in which it appears as an example of factorization of free monoids.
Actually, as pointed out to one of us by D. Knuth in 2004 the
reference [2] (the Chen, Fox, Lyndon 1958 \cite{CFL58}-L.A. Bokut,
Y. Q. Chen, Y. Li) it does not contain explicitly this statement."},
actually Sch\"utzenberger cited Chen-Fox-Lyndon 1958 paper by a
mistake (!).

So to the best of our knowledge the only origin of the ``Central
result" is Shirshov 1958 paper, at least for Sch\"utzenberger and
his school.

Some new linear bases of  $Lie(X)$ are given in the papers
\cite{BrKoSt02,Ci06,St08}.

In this paper we give another approach to definition of LS basis and
LS words following Shirshov's (Gr\"{o}bner-Shirshov bases) theory
for anti-commutative algebras \cite{Sh58a}. In our previous paper
\cite{BCL09} we gave the same kind of results for the Hall basis. To
be more precise, in that paper, we had found an anti-commutative
Gr\"{o}bner-Shirshov basis of a free anti-commutative
(non-associative) algebra $AC(X)$, such that the corresponding
irreducible monomials (not containing the maximal monomials of the
Gr\"{o}bner-Shirshov basis) are exactly the Hall monomials.

Here we prove the same kind of results for LS basis of $Lie(X)$.
Namely, we find a linear basis $N(X)$ of $AC(X)$ such that the
Composition-Diamond lemma is valid for the triple $(AC(X), N(X),
\succ),$ where $\succ$ is the deg-lex ordering on $N(X)$. So
Gr\"{o}bner-Shirshov bases theory is valid for the triple.

Then we find a Gr\"{o}bner-Shirshov basis of the ideal $J(X)$ of
$AC(X)$ generated by all Jacobians $J(N(X), N(X), N(X))$  on $N(X)$.
We do it in two steps. First we prove that the irreducible set of
monomials from $N(X)$ relative to $J(N(X), N(X), N(X))$  is exactly
the set of Lyndon-Shirshov monomials. It gives a new approach to LS
monomials. By the way, we also prove that for any monomial $(u)\in
N(X)$ its associative support $u$ has a form $v^m$, where $v$ is a
LS word. In some sense it gives a new  approach to LS words.

Our next step is the main result of the paper (see Theorem
\ref{t4.6}):

Let  $(AC(X), N(X), \succ, J(X))$ be as above. Then the set
$J(LS(X), LS(X), LS(X))$ of all Jacobians on LS monomials is a
Gr\"{o}bner-Shirshov basis of the ideal $J(X)$. What is more,
$Irr(J(LS(X), LS(X), LS(X)))$, the set of all monomials of $N(X)$
not containing maximal monomials of polynomials from  $J(LS(X),
LS(X), LS(X))$, is exactly the set of all LS monomials. Hence, the
later is a linear basis of the free Lie algebra $Lie(X)$.

For completeness we formulate Composition-Diamond lemma for the
triple $(AC(X), N(X),\\ \succ)$. The proof is essentially the same
as in \cite{Sh58} and \cite{BCL09} in which the difference is that
now $N(X)$ and the ordering are different from those papers.

\section{A basis of a free anti-commutative algebra AC(X)}
Let $X$ be a well-ordered set, $X^{*}$ the set of all associative
words, $X^{**}$ the set of all nonassociative words, $>_{lex}$,
$>_{dex-lex}$ the lexicographical ordering and the degree
lexicographical ordering on $X^{*}$ respectively.

Let $(u)$ be a nonassociative word on $X$ and denote $|(u)|$ the
degree of $(u)$, the number of the letters occur in $(u)$.

An ordering $\succ_{lex}$ on $X^{**}$ is inductively defined by:
$((u_1)(u_2))\succ_{lex}((v_1)(v_2))$ if and only if one of the following cases
holds

(a) $u_1u_2>_{lex}v_1v_2$,

(b) $u_1u_2=v_1v_2$ and $(u_1)\succ_{lex}(v_1)$,

(c) $u_1u_2=v_1v_2$, $(u_1)=(v_1)$ and $(u_2)\succ_{lex}(v_2).$

We define normal words $(u)\in X^{**}$ by induction on $|(u)|$:

(i) $x_i\in X$ is a normal word.

(ii) $(u)=((u_1)(u_2))$ is normal if both $(u_1)$ and $(u_2)$ are
normal and $(u_1)\succ_{lex}(u_2)$.

Denote $(u)$  by $[u]$, if $(u)$ is normal.\ \

Before we discuss some combinatorial properties of normal words on
$X$, we need to cite the definition of associative LS words and some
important properties of associative LS words.

An associative word $u$ is an associative LS word (ALSW for short)
if, for arbitrary nonempty $v$ and $w$ such that $u=vw$, we have
$u>_{lex}wv$.

We will use the following properties of ALSWs (see \cite{BC07}):

1. if $u = vw$ is an ALSW, where $v,w$ are not empty word, then
$u>_{lex}w$;

2. if $u, v$ are ALSWs and $u >_{lex} v$, then $uv$ is also an ALSW;

3. each nonempty word $u\in X^{*}$ has a unique decomposition:
$$
u=u_1u_2\ldots u_n,
$$

where $u_i$ is an ALSW and $u_1\leq_{lex}
u_2\leq_{lex}\ldots\leq_{lex} u_n$.

\begin{lemma}\label{l1} Let $u,v$ be ALSWs and $u>_{lex}v$. Then
$u^m>_{lex}v^n$ for any $m,n\geq 1$.
\end{lemma}
\noindent\textbf{Proof.}  If $u=x_1\cdots x_ty\cdots$ and
$v=x_1\cdots x_tz\cdots$ , where $y>_{lex}z$, then $u^m>_{lex}v^n$
clearly for any $m,n\geq 1$. Suppose $u$ is a proper prefix of $v$,
i.e., $v=uw$ for some nonempty word $w$. Then there exists an $l$
such that $v=u^lc$, $u$ is not a prefix of $c$. If $c$ is empty,
then $v=u^l$, which leads to a contradiction by the property 3. This
means that $c$ is not empty. Then, by the property 1,
$u>_{lex}v>_{lex}c$ and thus $u=x_1\cdots x_ty\cdots$ and
$c=x_1\cdots x_tz\cdots$, where $t\geq 0$ and $y>_{lex}z$. If this
is the case, $u^m>_{lex}v^n$ clearly.\ \ \ \ $\square$

\begin{lemma}\label{l2} Let $[u]$ be a normal word on $X$. Then $u=w^n$, $w$
is an ALSW and $n\geq 1$.
\end{lemma}
\noindent\textbf{Proof.}  If $[u]=x_i$, then this is a trivial case.
Suppose $[u]=[[u_1][u_2]]$. By induction on $|[u]|$, we have
$u_1=w_1^{m_1}$ and $u_2=w_2^{m_2}$, where $w_1,\ w_2$ are ALSWs and
$m_1, m_2\geq 1$. If $w_1>_{lex}w_2$, then, by the property 1 and
the property 2 of ALSWs, $u=w_1^{m_1}w_2^{m_2}$ is an ALSW. If
$w_1=w_2$, then $u=w_1^{m_1+m_2}$. If $w_1<_{lex}w_2$, then, by
Lemma \ref{l1}, $u=w_1^{m_1}<_{lex}v=w_2^{m_2}$ and this contradicts
$[u]$ is a normal word.\ \ \ \ $\square$

\begin{corollary}\label{c1}
If $[u]$ is normal and $[u]=[[u_1][u_2]]$, then
$u_1u_2\geq_{lex}u_2u_1$.
\end{corollary}

Let $k$ be a field, $N(X)$ the set of all normal words $[u]$ on $X$
and $AC(X)$ a $k$-vector space spanned by $N(X)$. Now define the
product of normal words by the following way: for any
 $[u],[v]\in N(X)$,
\begin{equation*}
\lbrack u][v]=\left\{
\begin{array}{r@{\quad}l}
[[u][v]], & \ \ \ \mbox{ if } \  [u]\succ_{lex}[v] \\
-[[v][u]], & \ \ \ \mbox{ if } \   [u]\prec_{lex}[v] \\
0,\text{ \ \ } & \ \ \ \mbox{ if } \   [u]=[v]%
\end{array}%
\right.
\end{equation*}

\noindent \textbf{Remark:} \ By definition, for any $(u)\in X^{**}$,
there exists a unique $[v]\in N(X)$ such that, in $AC(X)$, $(u)=\pm
[v]$ or 0. Now we denote $[v]$ by $\widetilde{(u)}$ if $(u)\neq0$.

Then we can get the following theorem by a straightforward proof.
\begin{theorem}$AC(X)$ is a free
anti-commutative algebra generated by $X$.
\end{theorem}

A well ordering $>$ on $N(X)$ is called monomial if it satisfies the
following condition:
$$
[u]>[v]\Rightarrow \widetilde{([u][w])}>\widetilde{([v][w])}
$$
for any $[w]\in N(X)$ such that $[w]\ne [u],\ [w]\ne [v]$.

If  $>$ is a monomial ordering on $N(X)$, then we have
$$
[u]>[v]\Rightarrow [a[u]b]>\widetilde{(a[v]b)},
$$
where $[a[u]b]$ is a normal word with subword $[u]$ and
$(a[v]b)=[a[u]b]|_{[u]\mapsto [v] }\neq 0$.

We define an ordering $\succ_{deg-lex}$ on $N(X)$ by:
$((u_1)(u_2))\succ_{deg-lex}((v_1)(v_2))$ if and only if one of the
following cases holds

(a) $u_1u_2>_{deg-lex}v_1v_2$,

(b) $u_1u_2=v_1v_2$ and $(u_1)\succ_{deg-lex}(v_1)$,

(c) $u_1u_2=v_1v_2$, $(u_1)=(v_1)$ and $(u_2)\succ_{deg-lex}(v_2).$

In this paper, unless another statement, we use the above ordering
$\succ_{deg-lex}$ on $N(X)$.

Then we have the following lemmas.

\begin{lemma}\label{l3}
Suppose that $[ u ]_{\mu}$ and $[ u ]_{\nu}$ are two different
bracketing on a word $u$ such that both $[u]_{\mu}$ and $[u]_{\nu}$
are normal words. If $[u]_{\mu}$ $\succ_{deg-lex}$ $[u]_{\nu}$, then
$[u]_{\mu}$ $\prec_{lex}$ $[u]_{\nu}$.
\end{lemma}
\noindent\textbf{Proof.} Let us consider the  basic case
$u=x_ix_jx_k, i>j>k$.  Then $[u]_{\mu}$ $\succ_{deg-lex}$
$[u]_{\nu}$ implies that $[u]_{\mu}=[[x_ix_j]x_k]$ and
$[u]_{\nu}=[x_i[x_jx_k]]$, and thus $[u]_{\mu}=[[x_ix_j]x_k]$
$\prec_{lex}$ $[u]_{\nu}=[x_i[x_jx_k]]$. Suppose that
$[u]_{\mu}=[[u_1]_{\mu_1}[u_2]_{\mu_2}]$ and
$[u]_{\nu}=[[u_1']_{\nu_1}[u_2']_{\nu_2}]$. If $u_1>_{deg-lex}u_1'$,
then $u_1'$ is a proper prefix of $u_1$ and obviously $[u]_{\mu}$
$\prec_{lex}$ $[u]_{\nu}$. If $u_1=u_1'$ and
$[u_1]_{\mu_1}$$\succ_{deg-lex}$ $[u_1']_{\nu_1}$, then, by
induction on $|u|$, $[u_1]_{\mu_1}$$\prec_{lex}$ $[u_1']_{\nu_1}$
and thus $[u]_{\mu}$ $\prec_{lex}$ $[u]_{\nu}$. If
$[u_1]_{\mu_1}=[u_1']_{\nu_1}$, then $u_2=u_2'$ and $[u_2]_{\mu_2}$
$\succ_{deg-lex}$ $[u_2']_{\nu_2}$, and hence we complete our proof
by induction on $|u|$.\ \ \ \ $\square$

\begin{lemma}\label{l4}
The ordering $\succ_{deg-lex}$ is a monomial well ordering on
$N(X)$.
\end{lemma}
\noindent\textbf{Proof.}  It is easy to check that the ordering
$\succ_{deg-lex}$ is a well ordering. Suppose that $[u],[v],[w]\in
N(X)$, $[u]\succ_{deg-lex} [v]$, $[w]\neq [u]$ and $[w]\neq [v]$. If
$|[u]|>|[v]|$, then $\widetilde{[u][w]}\succ_{deg-lex}
\widetilde{[v][w]}$ clearly. If $|[u]|=|[v]|$ and $u>_{lex}v$, then
we consider the following three cases $w\geq_{lex}u>_{lex}v$,
$u>_{lex}w>_{lex}v$ and $u>_{lex}v\geq_{lex}w$. For the cases
$w\geq_{lex}u>_{lex}v$ and $u>_{lex}v\geq_{lex}w$, it is obvious
that $\widetilde{[u][w]}\succ_{deg-lex} \widetilde{[v][w]}$. For the
case $u>_{lex}w>_{lex}v$,
$\widetilde{[u][w]}=[[u][w]]\succ_{deg-lex}
[[w][v]]=\widetilde{[v][w]}$ since by Corrollary \ref{c1} we have
$uw\geq_{lex}wu>_{lex}wv$. If $u=v>_{lex}w$ or $w>_{lex}u=v$,
$\widetilde{[u][w]}\succ_{deg-lex} \widetilde{[v][w]}$ clearly.
Suppose that $u=v=w$ and $[u]\succ_{dex-lex}[v]\succ_{dex-lex}[w]$.
Then, by Lemma \ref{l3}, $[u]\prec_{lex}[v]\prec_{lex}[w]$ and thus
$\widetilde{[u][w]}=[[w][u]]\succ_{deg-lex}
[[w][v]]=\widetilde{[v][w]}$. We can similarly check the other two
cases $[u]\succ_{dex-lex}[w]\succ_{dex-lex}[v]$ and
$[w]\succ_{dex-lex}[u]\succ_{dex-lex}[v]$ to complete our proof.\ \
\ \ $\square$

\section{Composition-Diamond lemma for $AC(X)$}

In this section, we formulate Composition-Diamond lemma for the free
anti-commutative algebra $AC(X)$.

Given a polynomial $f\in AC(X)$, it has the leading word
$[\overline{f}] \in N(X)$ according to the ordering
$\succ_{dex-lex}$ on $N(X)$, such that
$$
f=\sum_{[u]\in N}f([u])[u]=\alpha
[\overline{f}]+\sum{\alpha}_i[u_i],
$$
where $[\overline{f}]\succ_{deg-lex} [u_i], \ \alpha , \ {\alpha}_i,
\ f([u])\in k$. We call $[\overline{f}]$ the leading term of $f$.
$f$ is called monic if $\alpha=1$.

\begin{definition}
Let $S\subset AC(X)$ be a set of monic polynomials, $s\in S$ and $(u)\in X^{**}$. We define $S$-word $%
(u)_s$ by induction:
\begin{enumerate}
\item[(i)] $(s)_s=s$ is an $S$-word of $S$-length
1.
\item[(ii)] If $(u)_s$ is an $S$-word of $S$-length k and $(v)$
is a nonassociative word of length $l$, then
\begin{equation*}
(u)_s(v)\ and\ (v)(u)_s
\end{equation*}
are $S$-words of length $k+l$.
\end{enumerate}
\end{definition}

\begin{definition}
An $S$-word $(u)_s$ is called a normal $S$-word, if $(u)_{[\bar s]
}=(a[\bar s ]b)$ is a normal word. We denote $(u)_s$ by $[u]_s$ if
$(u)_s$ is a normal $S$-word. We also call the normal $S$-word
$[u]_s$ to be normal $s$-word. From Lemma \ref{l3} it follows that
$\overline{[u]_s}=[u]_{[\bar s] }$.
\end{definition}

Let $f,g$ be monic polynomials in $AC(X)$. Suppose that there exist
$a,b\in {X^*} $ such that $[\bar f] =[a[\bar g] b]$, where $[agb]$
is a normal $g$-word. Then we set $[w]=[\bar f ]$ and define the
composition of inclusion
\begin{equation*}
(f,g)_{[w]}=f-[agb].
\end{equation*}
We note that
$
(f,g)_{[w]}\in Id(f,g)\ \ and\ \
\overline{(f,g)_{[w]}}\prec_{deg-lex}[w],
$
where $Id(f,g)$ is the ideal of $AC(X)$ generated by $f,g$.

Transformation $f\mapsto f-[agb]$ is called the Elimination of
Leading Word (ELW) of $g$ in $f$.

Given a nonempty subset $S\subset AC(X)$, we shall say that the
composition $(f,g)_{[w]}$ is trivial modulo $(S,[w])$, if
\begin{equation*}
(f,g)_{[w]}=\sum\limits_i\alpha_i[a_is_ib_i],
\end{equation*}
where each $\alpha_i\in k,\ a_i,b_i\in X^*,\ s_i\in S,\ [a_is_ib_i]$
is normal S-word and $[a_i[\bar{s_i}]b_i]\prec_{deg-lex}{[w]}$. If
this is the case, then we write $(f,g)_{[w]}\equiv 0\ mod
(S,{[w]})$.

Let us note that if $(f,g)_{[w]}$ goes to $0$ by ELW's of $S$,
then $(f,g)_{[w]}\equiv 0$ mod$(S,[w])$. Indeed, using ELW's of
$S$, we have
$$
(f,g)_{[w]}\mapsto (f,g)_{[w]}-\alpha_1[a_1s_1b_1]=f_2 \mapsto
f_2-\alpha_2[a_2s_2b_2]\mapsto \cdots \mapsto 0.
$$
So, $(f,g)_{[w]}=\sum\limits_i\alpha_i[a_is_ib_i]$ , where
$[a_i[\bar{s_i}]b_i] \preceq_{deg-lex}
\overline{(f,g)_{[w]}}\prec_{deg-lex}{[w]}$.

In general, for $p,q\in AC(X)$, we write
$$
p\equiv q\quad mod(S,[w])
$$
which means that $p-q=\sum\alpha_i [a_i s_i b_i] $, where each
$\alpha_i\in k, \ a_i,b_i\in X^{*}, \ s_i\in S$ and $[a_i [\bar
{s_i}] b_i]\prec_{deg-lex}[w]$.

\begin{definition}
Let $S\subset AC(X)$ be a nonempty set of monic polynomials and the
ordering $\succ_{deg-lex}$ on $N(X)$ as before. Then the set $S$ is
called a Gr\"{o}bner-Shirshov basis (GSB for short) if any
composition $(f,g)_{[w]}$ with $f,g\in S$ is trivial modulo
$(S,{[w]})$, i.e., $(f,g)_{[w]}\equiv0$ $mod(S,{[w]})$.
\end{definition}

\begin{lemma}\label{3.8}
Let $S\subset AC(X)$ be a nonempty set of monic polynomials and
$Irr(S)=\{[u]\in N(X) |[u]\ne [a[\bar s] b]\ a,b\in X^*,\ s\in S
\mbox{ and } [as b] \mbox{ is a normal } S\mbox{-word}\}$. Then for
any $f\in AC(X)$,
\begin{equation*}
f=\sum\limits_{[u_i]\preceq_{deg-lex} [\bar f] }\alpha_i[u_i]+
\sum\limits_{[a_j[\overline{s_j}]b_j]\preceq_{deg-lex}[\bar
f]}\beta_j[a_js_jb_j],
\end{equation*}
where each $\alpha_i,\beta_j\in k, \ [u_i]\in Irr(S)$ and
$[a_js_jb_j]$ normal $S$-word.
\end{lemma}

By a similar proof to the Theorem 3.10 in \cite{BCL09}, we have the
following theorem. Here the difference is that now $N(X)$ and the
ordering are different from \cite{BCL09}. We omit the detail.

\begin{theorem}
Let $S\subset AC(X)$ be a nonempty set of monic polynomials, the
ordering $\succ_{deg-lex}$ on $N(X)$ as before and $Id(S)$ the ideal
of $AC(X)$ generated by $S$. Then the following statements are
equivalent:
\begin{enumerate}
\item [(i)] $S$ is a Gr\"{o}bner-Shirshov basis.

\item [(ii)] $f\in Id(S)\Rightarrow [\bar f] =[a[\bar s ]b]$ for some $s\in S\
and\ a,b\in X^*$, where $[as b]$ is a normal $S$-word.

\item [(iii)] $Irr(S)=\{[u]\in N(X) |[u]\ne [a[\bar s] b]\ a,b\in X^*,\ s\in S \mbox{ and }
[as b] \mbox{ is a normal } S\mbox{-word}\}$ is a linear basis of
the algebra $AC(X|S)=AC(X)/Id(S)$.
\end{enumerate}
\end{theorem}

\section{Gr\"{o}bner-Shirshov basis for a free Lie algebra}

In this section, we  give a new approach to non-associative
Lyndon-Shirshov words. We represent the free Lie algebra by the free
anti-commutative algebra and give a Gr\"{o}bner-Shirshov basis for
the free Lie algebra.

The proof of the following theorem is straightforward and we omit
the detail.

\begin{theorem}
Let $AC(X)$ be the free anti-commutative algebra and let
\begin{equation*}
S=\{([u][v])[w]-([u][w])[v]-[u]([v][w]) \ | \ [u],[v],[w]\in N(X)
\mbox{\ and \ } [u]\succ_{lex}[v]\succ_{lex}[w]\}.
\end{equation*}
Then the algebra $AC(X|S)$ is the free Lie algebra generated by $X$.
\ \ \ \ $\square$
\end{theorem}

We now cite the definition of non-associative Lyndon-Shirshov words
(NLSWs) given by induction on length:
\begin{enumerate}
\item[1)]$x$ is a NLSW  for any $x\in X$,

\item[2)] a non-associative word $((v)(w))$ is called a  NLSW if
\begin{enumerate}

\item[(a)] both $(v)$ and $(w)$ are NLSWs, and $v>_{lex}w$,

\item[(b)] if $(v)=((v_1)(v_2))$, then $v_2\leq_{lex} w$.
\end{enumerate}
\end{enumerate}

Note that each NLSW $(u)$ is a normal word. We denote $(u)$ by
$[[u]]$ if $(u)$ is a NLSW.

There are some known properties of NLSW (see \cite{BC07}):

(1) for any NLSW $[[u]]$, $u$ is an ALSW;

(2) for each ALSW $u$, there exists a unique bracketing $(u)$ on $u$
 such that $(u)$ is a NLSW. It follows that
$[[u]]\succ_{lex}[[v]]$ iff $u>_{lex}v$.

In general, we don't know what is the leading term of a polynomial
in $S$. However, for some special polynomials, we know exactly the
leading terms of them, which can help us to study the ideal
generated by  $S$.

\begin{lemma}\label{4.1a}
Let $[[u]],[[v]],[[w]]$ be NLSWs and
$[[u]]\succ_{lex}[[v]]\succ_{lex}[[w]]$. Denote $f_{uvw}$ the
polynomial $([[u]][[v]])[[w]]-([[u]][[w]])[[v]]-[[u]]([[v]][[w]])$.
Then $\overline{f_{uvw}}=([[u]][[v]])[[w]]$.
\end{lemma}
\noindent\textbf{Proof.} Following from
$[[u]]\succ_{lex}[[v]]\succ_{lex}[[w]]$, we get $u>_{lex}v>_{lex}w$.
Then $\widetilde{([[u]][[v]])[[w]]}=([[u]][[v]])[[w]]$ $(\mbox{ for
} u>_{lex}uv>_{lex}v>_{lex}w)$ and
$\widetilde{[[u]]([[v]][[w]])}=[[u]]([[v]][[w]])$ $(\mbox{ for }
u>_{lex}v>_{lex}vw>_{lex}w)$. $(\widetilde{([[u]][[w]])[[v]]}$ may
be $([[u]][[w]])[[v]]$ or $[[v]]([[u]][[w]])$, however, we have
$uvw>uwv,vuw$ and thus $\overline{f_{uvw}}=([[u]][[v]])[[w]]$.\ \ \
\ $\square$

The following theorem gives us a new approach to NLSWs.

\begin{theorem}\label{T4.1}
Let $T$ be the set consisting of all NLSWs. Then
\begin{equation*}
Irr(S)=T.
\end{equation*}
\end{theorem}

\noindent\textbf{Proof.} Suppose $[u]\in Irr(S)$. We will show that
$[u]$ is a NLSW by induction on $|[u]|=n$. If $n=1$, then $[u]=x\in
X$ which is already a NLSW. Let $n>1$ and $[u]=[[v][w]]$. By
induction, we have that $ [v],[w] $ are NLSWs. If $|v|=1$, then
$[u]$ is a NLSW. If $|[v]|>1\ and \ [v]=[[v_1][v_2]]$, then
$v_2\leq_{lex} w$ because of Lemma \ref{4.1a}  and $[u]\in Irr(S)$.
So, $[u]$ is a NLSW.

It's clear that every NLSW is in $Irr(S)$ since every subword of
NLSW is also a NLSW.\ \ \ \ $\square$

In order to prove that $S$ is a GSB in $AC(X)$, we consider the
subset of $S$:
\begin{eqnarray*}
 S_0&=&\{([[u]][[v]])[[w]]-([[u]][[w]])[[v]]-[[u]]([[v]][[w]]) \ |
 \\
&&[[u]]\succ_{lex}[[v]]\succ_{lex}[[w]],\  and  \ [[u]],[[v]],[[w]]\
are\ NLSWs\}.
\end{eqnarray*}

In the following, we prove that $S_0$ generates the same ideal as
$S$, and $S_0$ is a GSB in $AC(X)$, which implies that $S$ is also a
GSB.

According to Lemma \ref{3.8} and $Irr(S_0)=T$ (similar proof to Theorem \ref{T4.1}), we get the following lemma.
\begin{lemma}\label{4.2}
In $AC(X)$, any normal word $[u]$ has the following presentation:
$$
[u]=\sum\limits_{i}\alpha_i[[u_i]]+\sum\limits_{j}\beta_j[u_j]_{s_j},
$$
where $\alpha_i, \beta_j\in k$, $[[u_i]]$ are NLSWs, $[u_j]_{s_j}$
normal $S_0$-words, $s_j\in S_0, \ [[u_i]], [u_j]_{[\overline
{s_j}]}\preceq_{deg-lex} [u]$. Moreover, each $[[u_i]]$ has the same
degree as $[u]$.
\end{lemma}

\begin{lemma}
Suppose $S$ and $S_0$ are sets defined as before. Then, in $AC(X)$,
we have
\begin{equation*}
Id(S)=Id(S_0).
\end{equation*}

\end{lemma}
{\bf Proof.} Since $S_0$ is a subset of $S$, it suffices to prove
that $AC(X|S_0)$ is a Lie algebra. We need only to prove that, in
$AC(X|S_0)$,
\begin{equation*}
([u][v])[w]-([u][w])[v]-[u]([v][w])=0,
\end{equation*}
where $[u],[v],[w]\in N(X)$. By Lemma \ref{4.2}, it suffices to
prove
\begin{equation*}
([[u]][[v]])[[w]]-([[u]][[w]])[[v]]-[[u]]([[v]][[w]])=0,
\end{equation*}
where  $[[u]]\succ_{lex}[[v]]\succ_{lex}[[w]]$. This is trivial by
the definition of $S_0$. \ \ \ \ $\square$

\begin{theorem}\label{t4.6}
Let
\begin{center}
 $S_0=\{([[u]][[v]])[[w]]-([[u]][[w]])[[v]]-[[u]]([[v]][[w]]) \ |
~[[u]]\succ_{lex}[[v]]\succ_{lex}[[w]],\  and  \ [[u]],[[v]],[[w]]\
\mbox{ are nonassociative LS words}\ \}.$
\end{center}
 Then the set $S_0$ is a Gr\"obner-Shirshov basis in
$AC(X)$.
\end{theorem}
{\bf Proof.} To simplify notations, we write $\hat{u}$ for $[[u]]$ and
$\hat{u_1}\hat{u_2}\cdots \hat{u_n}$ for
$((((\hat{u_1})\hat{u_2})\cdots)\hat{u_n})$.\\
Let $
f_{uvw}=\hat{u}\hat{v}\hat{w}-\hat{u}\hat{w}\hat{v}-\hat{u}(\hat{v}\hat{w}),
$ where $\hat{u},\hat{v},\hat{w} $ are NLSWs and
$u>_{lex}v>_{lex}w$. We know
$\overline{f_{uvw}}=\hat{u}\hat{v}\hat{w}$ from Lemma \ref{4.1a}.

Suppose $\overline{f_{u_1v_1w_1}}$ is a subword of
$\overline{f_{uvw}}$. Since $\hat{u},\hat{v},\hat{w}$ are NLSWs, we
have
$\hat{u_1}\hat{v_1}\hat{w_1}=\hat{u}\hat{v},\hat{u}=\hat{u_1}\hat{v_1}\
\mbox { and } \ \hat{v}=\hat{w_1}$. We will prove that the
composition
\begin{equation*}
(f_{uvw},f_{u_1v_1w_1})_{\hat{u}\hat{v}\hat{w}}
\end{equation*}
is trivial modulo $(S_0, \hat{u}\hat{v}\hat{w})$. We note that
$u_1>_{lex}v_1>_{lex}w_1=v>_{lex}w$.

Firstly, we  prove that the following statements hold
mod$(S_0,\hat{u}\hat{v}\hat{w})$:
\begin{enumerate}

\item[1)]$\hat{u_1}\hat{v}\hat{v_1}\hat{w}-\hat{u_1}\hat{v}\hat{w}\hat{v_1}-\hat{u_1}\hat{v}(\hat{v_1}\hat{w})\equiv 0$.

\item[2)]$\hat{u_1}(\hat{v_1}\hat{v})\hat{w}-\hat{u_1}\hat{w}(\hat{v_1}\hat{v})-\hat{u_1}(\hat{v_1}\hat{v}\hat{w})\equiv 0$.

\item[3)]$\hat{u_1}\hat{w}\hat{v_1}\hat{v}-\hat{u_1}\hat{w}\hat{v}\hat{v_1}-\hat{u_1}\hat{w}(\hat{v_1}\hat{v})\equiv 0$.

\item[4)]$\hat{u_1}(\hat{v_1}\hat{w})\hat{v}-\hat{u_1}\hat{v}(\hat{v_1}\hat{w})-\hat{u_1}(\hat{v_1}\hat{w}\hat{v})\equiv 0$.

\item[5)]$\hat{u_1}(\hat{v}\hat{w})\hat{v_1}-\hat{u_1}(\hat{v}\hat{w}\hat{v_1})-\hat{u_1}\hat{v_1}(\hat{v}\hat{w})\equiv 0$.

\item[6)]$\hat{u_1}\hat{v}\hat{w}\hat{v_1}-\hat{u_1}\hat{w}\hat{v}\hat{v_1}-\hat{u_1}(\hat{v}\hat{w})\hat{v_1}\equiv 0$.

\item[7)]$\hat{u_1}\hat{v_1}\hat{w}\hat{v}-\hat{u_1}\hat{w}\hat{v_1}\hat{v}-\hat{u_1}(\hat{v_1}\hat{w})\hat{v}\equiv 0 $.

\item[8)]$\hat{u_1}(\hat{v_1}\hat{v}\hat{w})-\hat{u_1}(\hat{v_1}\hat{w}\hat{v})-\hat{u_1}(\hat{v_1}(\hat{v}\hat{w}))\equiv 0$.

\end{enumerate}

We only prove $1)$. $2)$--$5)$ are similarly proved to $1)$ and
$6)$--$8)$ follow from ELW's of $S_0$.

By ELW's of $S_0$, we may assume, without loss of generality, that
$\hat{u_1}\hat{v}$ is a NLSW. If $u_1v>_{lex}v_1$, then
$\hat{u_1}\hat{v}\hat{v_1}\hat{w}-\hat{u_1}\hat{v}\hat{w}\hat{v_1}-\hat{u_1}\hat{v}(\hat{v_1}\hat{w})
=f_{(u_1v)v_1w}\equiv 0$. If $u_1v=v_1$, then
$\hat{u_1}\hat{v}\hat{v_1}\hat{w}-\hat{u_1}\hat{v}\hat{w}\hat{v_1}-\hat{u_1}\hat{v}(\hat{v_1}\hat{w})
= 0$. If $u_1v<_{lex}v_1$, then
$\hat{u_1}\hat{v}\hat{v_1}\hat{w}-\hat{u_1}\hat{v}\hat{w}\hat{v_1}-\hat{u_1}\hat{v}(\hat{v_1}\hat{w})
=-f_{v_1(u_1v)w}\equiv 0$.

Secondly, we have
\begin{align*}
(f_{uvw},f_{u_1v_1w_1})_{\hat{u}\hat{v}\hat{w}}&=f_{uvw}-(f_{u_1v_1w_1})\hat{w } \\
&=\hat{u_1}\hat{v}\hat{v_1}\hat{w}+\hat{u_1}(\hat{v_1}\hat{v})\hat{w}-\hat{u_1}\hat{v_1}\hat{w}\hat{v}-\hat{u_1}\hat{v_1}(\hat{v}\hat{w}).
\end{align*}
Let
\begin{equation*}
A=\hat{u_1}\hat{v}\hat{v_1}\hat{w}+\hat{u_1}(\hat{v_1}\hat{v})\hat{w} \ \ and \ \ B=-\hat{u_1}\hat{v_1}\hat{w}\hat{v}-\hat{u_1}\hat{v_1}(\hat{v}\hat{w}).
\end{equation*}

Then, by $1)$--$8)$, we have
\begin{align*}
A & \equiv \hat{u_1}\hat{v}\hat{w}\hat{v_1}+\hat{u_1}\hat{v}(\hat{v_1}\hat{w})+\hat{u_1}\hat{w}(\hat{v_1}\hat{v})+\hat{u_1}(\hat{v_1}\hat{v}\hat{w}) \ \ \ \\
&\equiv
\hat{u_1}\hat{w}\hat{v}\hat{v_1}+\hat{u_1}(\hat{v}\hat{w})\hat{v_1}+\hat{u_1}\hat{v}(\hat{v_1}\hat{w})+\hat{u_1}\hat{w}(\hat{v_1}\hat{v})+\hat{u_1}(\hat{v_1}\hat{w}\hat{v})+\hat{u_1}(\hat{v_1}(\hat{v}\hat{w})) \
\end{align*}
and
\begin{align*}
-B&=\hat{u_1}\hat{v_1}\hat{w}\hat{v}+(\hat{u_1}\hat{v_1})(\hat{v}\hat{w}) \\
&\equiv \hat{u_1}\hat{w}\hat{v_1}\hat{v}+\hat{u_1}(\hat{v_1}\hat{w})\hat{v}+\hat{u_1}\hat{v_1}(\hat{v}\hat{w}) \ \\
&\equiv \hat{u_1}\hat{w}\hat{v}\hat{v_1}+\hat{u_1}\hat{w}(\hat{v_1}\hat{v})+\hat{u_1}\hat{v}(\hat{v_1}\hat{w})+\hat{u_1}(\hat{v_1}\hat{w}\hat{v})+\hat{u_1}\hat{v_1}(\hat{v}\hat{w}). \ \
\
\end{align*}
So,
\begin{equation*}
(f_{uvw},f_{u_1v_1w_1})_{\hat{u}\hat{v}\hat{w}}=A+B\equiv
\hat{u_1}(\hat{v}\hat{w})\hat{v_1}+\hat{u_1}(\hat{v_1}(\hat{v}\hat{w}))-\hat{u_1}\hat{v_1}(\hat{v}\hat{w})\equiv0
\ \ mod(S_0, \hat{u}\hat{v}\hat{w}).
\end{equation*}

This completes our proof. \ \ \ \ $\square$

\end{document}